# Several new domain-type and boundary-type numerical discretization schemes with radial basis function

W. Chen

*Department of Informatics, University of Oslo, P.O.Box 1080, Blindern, 0316 Oslo, Norway*

**Abstract**

This paper is concerned with a few novel RBF-based numerical schemes discretizing partial differential equations. For boundary-type methods, we derive the indirect and direct symmetric boundary knot methods (BKM). The resulting interpolation matrix of both is always symmetric irrespective of boundary geometry and conditions. In particular, the direct BKM applies the practical physical variables rather than expansion coefficients and becomes very competitive to the boundary element method. On the other hand, based on the multiple reciprocity principle, we invent the RBF-based boundary particle method (BPM) for general inhomogeneous problems without a need using inner nodes. The direct and symmetric BPM schemes are also developed.

For domain-type RBF discretization schemes, by using the Green integral we develop a new Hermite RBF scheme called as the modified Kansa method (MKM), which differs from the symmetric Hermite RBF scheme in that the MKM discretizes both governing equation and boundary conditions on the same boundary nodes. The local spline version of the MKM is named as the finite knot method (FKM). Both MKM and FKM significantly reduce calculation errors at nodes adjacent to boundary. In addition, the nonsingular high-order fundamental or general solution is strongly recommended as the RBF in the domain-type methods and dual reciprocity method approximation of particular solution relating to the BKM.

It is stressed that all the above discretization methods of boundary-type and domain-type are symmetric, meshless, and integration-free. The spline-based schemes will produce desirable symmetric sparse banded interpolation matrix. In appendix, we present a Hermite scheme to eliminate edge effect on the RBF geometric modeling and imaging.

## 1. Introduction

In recent years there is a strong interest in developing the meshfree numerical schemes due to the fact that in the traditional finite element and boundary element methods, the mesh generation of high-dimensional problems costs a huge amount of computer resources [1,2]. Almost all meshfree methods available today require using the moving least square (MLS) technique [3]. Exceptionally, the numerical schemes based on the radial basis function (RBF) do not need to use the MLS at all and are inherently meshfree. For some new advances on the RBF see Buhmann's excellent survey [4]. The RBF-type methods have physical backgrounds of field and potential theory [5] and are justified mathematically by integral equation theory, that is, the RBF techniques have underlying relationship with the Green integral and translation invariance kernel function of integral equation [6].

Since early 90's, some efficient RBF schemes for numerical PDEs have been presented. Among them famous are the Kansa method [7], Hermite domain-type RBF collocation method [8] and the method of fundamental solution (MFS) [9]. Most recently, the present author [5,6] developed the boundary knot method (BKM) as a competitive alternative to the MFS. In fact, earlier than Kansa's pioneer work [7], Nardini and Brebbia [10] in early 80's, without knowing the RBF terminology and existent developments, applied the RBF and dual reciprocity method to effectively eliminate domain integral in context of the boundary element technique. This original work gives rise to currently popular dual reciprocity BEM (DRBEM). All of the above RBF methods perform excellently in some given numerical experiments.

The Kansa's method [7,11] is the very first domain-type RBF collocation scheme with easy-to-use merit. However, the method loses the symmetric interpolation matrix due to the mixed boundary conditions and/or non-self-adjoint operator. The Hermite RBF collocation method [8] overcomes the unsymmetrical drawbacks in the Kansa's method. Like the Kansa's method, however, the method still suffers relatively low accuracy in boundary-adjacent region. The MFS is a simple and efficient boundary-type RBF scheme. The shortcoming is that the controversial artificial boundary outside physical domain impedes its practical applications [12]. The BKM [5,6] surpasses the MFS in that it employs the non-singular general solution instead of the MFS's singular fundamental solution. Therefore, there is no longer need to use the arbitrary fictitious boundary. It is found that the BKM can produce very accurate solution for complicated geometry problems with relatively fewer nodes [13]. On the other hand, it is noted that the MFS could not yield symmetric interpolation matrix, while the BKM, just like the Kansa method, loses symmetric merit when involving the mixed boundary conditions and/or non-self-adjoin operator. In addition, it is especially worth pointing out that all above RBF schemes are indirect and global. The indirect methods mean that the expansion coefficients rather than physical variables are used as the basic variable in the whole computing procedure, while the global interpolation causes the ill-conditioning interpolation matrix.

The purpose of this paper is to introduce a few new RBF discretization schemes of boundary and domain types, all of which produce the symmetric interpolation matrix no matter what kind of boundary conditions and geometry. The symmetricity guarantees the reliability of these methods with some conditions. In particular, we present the direct and indirect RBF schemes and their local spline versions.

The rest of this paper is structured as follows. In section 2, a few boundary-type RBF schemes are presented. Firstly, we establish the symmetric Hermite BKM, and then the direct BKM is introduced which uses the physical variables instead of expansion coefficients. The method has great potential to challenge the BEM as a typical boundary-type numerical technique. We also strongly recommend the nonsingular higher-order fundamental or general solution as the RBF for approximate expression of particular solution. Secondly, we derive the boundary particle method (BPM) by using the multiple reciprocity principle and RBF. The symmetric direct BPM is also proposed. It is noted that unlike the BKM, the BPM does not require the interior nodes to improve accuracy.

Section 3 is concerned with the domain-type RBF schemes. By using Green integral, we present the symmetric modified Kansa's method (MKM) which collocates both the governing and boundary equations at the same boundary nodes. This method significantly improves the solution accuracy at nodes neighboring boundary. Furthermore, we develop the spline version of the MKM called as the finite knot method (FKM), which produces the sparse symmetric interpolation matrix.

Finally, section 5 concludes some remarks on these novel RBF schemes. We proposed the spline version of the BKM and BPM to produce the sparse banded interpolation matrix. Rapid calculation and wavelets multiscale analysis of the RBF collocation are also briefly discussed. In appendix, we give a Hermite RBF interpolation scheme to eliminate edge effect of geometric reconstruction and imaging.

## 2. New boundary-type RBF schemes

### 2.1. Boundary knot methods

The present author recently introduced RBF-based and integration-free boundary knot method [5,6]. Some numerical experiments demonstrated the BKM performed very well for 2D and 3D Helmholtz, modified Helmholtz and convection-diffusion problems with very complicated geometry [13]. The aim of this study is to derive the symmetric Hermite BKM and direct BKM.

To clearly illustrate our idea, consider the following example without loss of generality

$$L\{u\} = f(x), x \in \Omega, \quad (1)$$

$$u(x) = D(x), x \subset S_u, \quad (2a)$$

$$\frac{\partial u(x)}{\partial n} = N(x), x \subset S_T, \quad (2b)$$

where $x$ means multi-dimensional independent variable, and $n$ is the unit outward normal. The solution of Eq. (1) can be expressed as

$$u = v + u_p, \quad (3)$$

where $v$ and $u_p$ are the general and particular solutions, respectively. The latter satisfies

$$L\{u_p\} = f(x) \quad (4)$$

but does not necessarily satisfy boundary conditions.

To evaluate the particular solution, the inhomogeneous term is approximated first by

$$f(x) \cong \sum_{j=1}^{N+L} \alpha_j \varphi(r_j), \quad (5)$$

where $\alpha_j$ are the unknown coefficients. $N$ and $L$ are respectively the numbers of knots on the domain and boundary. The use of interior points here is usually necessary to guarantee the accuracy and convergence of the BKM solution.

$r_j = \|x - x_j\|$ represents the Euclidean distance norm, and $\varphi$ is the radial basis function.

By forcing approximation representation (5) to exactly satisfy Eq. (4) at all nodes, we can uniquely determine

$$\alpha = A_\varphi^{-1} \{f(x_i)\}, \quad (6)$$

where $A_\varphi$ is nonsingular RBF interpolation matrix. Finally, we can get particular solutions at any point by summing localized approximate particular solutions

$$u_p = \sum_{j=1}^{N+L} \alpha_j \phi(\|x - x_j\|). \quad (7)$$

Substituting Eq. (6) into Eq. (7) yields

$$u_p = \Phi A_\varphi^{-1} \{f(x_i)\}, \quad (8)$$

where $\Phi$ is a known matrix comprised of $\phi(r_{ij})$. In most practical applications, the approximate particular solution $\phi$ is determined beforehand, and then we evaluate the corresponding $\varphi$ through substituting $\phi$ into differential operator $L$. Ref. [14] favors the thin plate spline (TPS) for the particular solution approximation. However, the TPS is only applicable to the 2D Laplace operator rather than to general cases due to its implicit linkage with the second order fundamental solution of 2D Laplace operator. Generally speaking, the non-singular high-order fundamental solution in the multiple reciprocity method [15] is strongly recommended as the radial basis function $\phi$ in both the dual reciprocity BEM (DRBEM) and the BKM.

On the other hand, the homogeneous solution $v$ has to satisfy

$$\nabla^2 v = 0, \quad (9)$$

$$v(x) = D(x) - u_p, \quad (10a)$$

$$\frac{\partial v(x)}{\partial n} = N(x) - \frac{\partial u_p(x)}{\partial n}. \quad (10b)$$

Unlike the dual reciprocity BEM [10,16] and MFS [9] using the singular fundamental solution, the BKM [5,6] approximates $v$ by means of nonsingular general solution, namely,

$$v(x) = \sum_{k=1}^{L} \lambda_k u^{\#}(r_k), \quad (11)$$

where $r_k = \|x - x_k\|$ and $k$ is index of source points on boundary; $u^{\#}$ is the nonsingular general solution of operator $L$. The non-singular general solutions of some frequently-used operators can be found in Chen and Tanaka [5,6]. $\lambda_k$ are the desired coefficients. Collocating Eqs. (10a,b) at all boundary and interior knots in terms of representation (11), we have

$$\sum_{k=1}^{L} \lambda_k u^{\#}(r_{ik}) = D(x_i) - u_p(x_i), \quad (12)$$

$$\sum_{k=1}^{L} \lambda_k \frac{\partial u^{\#}(r_{jk})}{\partial n} = N(x_j) - \frac{\partial u_p(x_j)}{\partial n}, \quad (13)$$

$$\sum_{k=1}^{L} \lambda_k u^{\#}(r_{lk}) = u_l - u_p(x_l), \quad l = 1,\ldots,N, \quad (14)$$

where $i$, $j$, and $l$ indicate response knots respectively located on boundary $S_u$, $S_\Gamma$, and domain $\Omega$. Substituting approximate particular solution (8) into Eqs. (12), (13) and (14), we can solve the above simultaneous algebraic equations. After this, we can employ the obtained expansion coefficients $\lambda_k$ and inner knot solutions $u_l$ to calculate the BKM solution at any knot.

### 2.1.1. Symmetric indirect BKM

As is pointed out in [5], the BKM interpolation matrix is symmetric only if one kind of boundary condition is involved in given problems. Otherwise, the symmetricity is destroyed as in the Kansa's method. By analogy with the symmetric Hermite RBF collocation method presented by Fasshauer [8], we modify the BKM approximate expression (11) of homogeneous solution $v$ as

$$v(x) = \sum_{s=1}^{L_d} a_s u^{\#}(r_s) - \sum_{t=1}^{L_N} b_t \frac{\partial u^{\#}(r_t)}{\partial n}, \quad (15)$$

where $n$ is the unit outward normal as in boundary condition (2b), and $L_d$ and $L_n$ are respectively the numbers of knots at the Dirichlet and Neumann boundary surfaces. The minus sign associated with the second term is due to the fact that the Neumann condition of the first order derivative is not self-adjoint. In terms of expression (15), the collocation analogue equations (10a,b) are rewritten as

$$\sum_{s=1}^{L_d} a_s u^{\#}(r_{is}) - \sum_{t=1}^{L_N} b_t \frac{\partial u^{\#}(r_{it})}{\partial n} = D(x_i) - u_p(x_i), \quad (16)$$

$$\sum_{s=1}^{L_d} a_s \frac{\partial u^{\#}(r_{js})}{\partial n} - \sum_{t=1}^{L_N} b_t \frac{\partial^2 u^{\#}(r_{jt})}{\partial n^2} = N(x_j) - \frac{\partial u_p(x_j)}{\partial n}, \quad (17)$$

$$\sum_{s=1}^{L_d} a_s u^{\#}(r_{ls}) - \sum_{t=1}^{L_N} b_t \frac{\partial u^{\#}(r_{lt})}{\partial n} = u_l - u_p(x_l), \quad (18)$$

The solution of the above simultaneous equations can be decomposed into two steps. The first is to evaluate the unknown boundary expansion coefficients $a_s$ and $b_t$ by using Eqs. (16) and (17), and then the interior node solution $u_l$ is calculated by Eq. (18). The system matrix of Eqs. (16) and (17)

$$A = \begin{bmatrix} u^{\#}(r_{is}) & -\dfrac{\partial u^{\#}(r_{it})}{\partial n} \\ \dfrac{\partial u^{\#}(r_{js})}{\partial n} & -\dfrac{\partial^2 u^{\#}(r_{jt})}{\partial n^2} \end{bmatrix} \quad (19)$$

is symmetric if operator $L\{\}$ is self-adjoint. Here we can see that minus sign keep the two minor-diagonal entries the same value. Note that $i$, $s$, $j$ and $t$ are reciprocal indices of Dirichlet and Neumann boundary nodes. It is stressed here that the MFS could not produce the symmetric interpolation matrix in any way.

It is worth pointing out that for more complicated mixed differential boundary conditions, we still can construct the corresponding Hermite BKM boundary expression for producing the symmetric boundary interpolation matrix. Therefore, the presented symmetric BKM is general enough to handle a broad range of practical problems.

### 2.1.2. Symmetric direct BKM

It is well known that the direct BEM is more popular than the older indirect BEM due to its strengths in some problems [17]. We note that the above BKM use the expansion coefficient rather than the direct physical variable in the approximation of boundary value. Therefore, such BKM is called as the indirect BKM.

Still considering the problem described by Eqs. (1) and (2a,b), the Neumann condition $N(x)$ at $x \subset S_u$ and Dirichlet condition $D(x)$ at $x \subset S_T$ are unknown in contrast to the prescribed boundary condition (2a,b). To simplify the presentation, we use the $D_u$, $N_u$ and $D_\Gamma$, $N_\Gamma$ respectively represent the Dirichlet and Neumann values at $x \subset S_u$ and $x \subset S_T$. In terms of Eqs. (15), (16) and (17), we have

$$\begin{bmatrix} u^\#(r_{is}) & -\dfrac{\partial u^\#(r_{it})}{\partial n} \\ \dfrac{\partial u^\#(r_{js})}{\partial n} & -\dfrac{\partial^2 u^\#(r_{jt})}{\partial n^2} \end{bmatrix} \begin{Bmatrix} a_s \\ b_t \end{Bmatrix} = \begin{Bmatrix} D_u \\ N_\Gamma \end{Bmatrix}, \quad (20)$$

$$\begin{bmatrix} \dfrac{\partial u^\#(r_{is})}{\partial n} & -\dfrac{\partial^2 u^\#(r_{it})}{\partial n^2} \\ u^\#(r_{js}) & -\dfrac{\partial u^\#(r_{jt})}{\partial n} \end{bmatrix} \begin{Bmatrix} a_s \\ b_t \end{Bmatrix} = \begin{Bmatrix} N_u \\ D_\Gamma \end{Bmatrix}. \quad (21)$$

In matrix form, Eqs. (20) and (21) are restated as

$$A \begin{Bmatrix} a_s \\ b_t \end{Bmatrix} = \begin{Bmatrix} D_u \\ N_\Gamma \end{Bmatrix}, \quad (22)$$

$$B \begin{Bmatrix} a_s \\ b_t \end{Bmatrix} = \begin{Bmatrix} N_u \\ D_\Gamma \end{Bmatrix}. \quad (23)$$

So we have

$$\begin{Bmatrix} N_u \\ D_\Gamma \end{Bmatrix} = BA^{-1} \begin{Bmatrix} D_u \\ N_\Gamma \end{Bmatrix}, \quad (24)$$

where the left-hand boundary values are unknown, while the right-hand boundary values are prescribed in Eqs. (1,2). In particular, it is noted that matrix $A$ are symmetric and no indirect expansion coefficients are involved in Eq. (24). Therefore, the above solution procedure is a symmetric direct BKM strategy. It is more straightforward to derive unsymmetrical direct BKM scheme.

Just like the comparisons between the direct BEM and indirect BEM [17], the direct BKM has the advantages over the indirect BKM in that it is more feasible to some problems with sharp corners since the fictitious expansion coefficients may tend to infinity as the nodes increase, even if the physical quantity remain well behaved. In addition, the discretization systems of the direct method is more efficiently solved by iterative techniques than that of the indirect method due to the fact that the former is easier to find nice initial guess solutions in the iterative algorithms [18].

On the other hand, it is worth stressing that if we employ the Fasshaure's Hermite scheme [8] instead of approximate expression (5) of particular solution, RBF interpolation matrix $A_\varphi$ in Eqs. (6) and (8) will also be symmetric. The symmetricity guarantees the solvability of the BKM methodology with some conditions.

**2.2. Boundary particle methods**

As in the dual reciprocity BEM, the interior nodes are usually necessary in the BKM. A rival to the DRBEM is the multiple reciprocity BEM often called shortly as the MRM [19], which applies the multiple reciprocity principle to circumvent the domain integral without using any inner nodes. The shortcoming in the MRM is uneasily applied to nonlinear problems and higher computing effort compared with the DRBEM. In this section, we will develop a new boundary-only RBF scheme based on the multiple reciprocity principle.

The MRM assumes that the particular solution in Eq. (3) can be approximated by higher-order homogeneous solution, namely,

$$u = v^0 + u_p^0 = v^0 + \sum_{m=1}^{\infty} v^m, \quad (25)$$

where superscript $m$ is the order index of homogeneous solution. $v^0$ and $u_p^0$ are equivalent to homogeneous solution $v$ and particular solution $u_p$ in Eq. (3). Through an incremental differentiation operation via operator $L\{\}$, we have successively higher order differential equations:

$$\begin{cases} v^0(x_i) = D(x_i) - u_p^0(x_i) \\ \dfrac{\partial v^0(x_j)}{\partial n} = N(x_j) - \dfrac{\partial u_p^0(x_j)}{\partial n} \end{cases}, \quad (26a)$$

$$L^0\{v^1(x)\} = f(x) - L^0\{u_p^1(x)\}, \quad (26b)$$

$$L^{n-1}\{v^n(x)\} = L^{n-2}\{f(x)\} - L^{n-1}\{u_p^n(x)\},$$
$$n=2,3,\ldots \quad (26c)$$

where $L^n\{\}$ denotes the $n$-th order of operator $L\{\}$, say $L^1\{\}=LL^0\{\}$, and $L^0\{\}$ equals $L\{\}$. For example, the $n$-th order of Laplace operator is $\nabla^{2(n+1)}$. $u_p^n$ is the $n$-th order of particular solution defined as

$$u_p^n = \sum_{m=n+1}^{\infty} v^m. \quad (27)$$

The homogeneous solution $v^m$ is approximated by

$$v^m(x) = \sum_{k=1}^{L} \lambda_k^m u_m^\#(r_k), \quad (28)$$

where $L$ is the number of boundary nodes, and $u_m^\#$ is the corresponding $m$-th order fundamental solution or nonsingular general solution satisfying

$$L^m\{u_m^\#\} = 0. \quad (29)$$

Collocating boundary equations (26a,b,c) only on boundary nodes, we have the boundary discretization equations

$$\left.\begin{aligned} \sum_{k=1}^{L} \lambda_k^0 u^0(r_{ik}) &= D(x_i) - u_p^0(x_i) \\ \sum_{k=1}^{L} \lambda_k^0 \dfrac{\partial u^0(r_{jk})}{\partial n} &= N(x_j) - \dfrac{\partial u_p^0(x_j)}{\partial n} \end{aligned}\right\} = b^0, \quad (30a)$$

$$\sum_{k=1}^{L} \lambda_k^1 L^0\{u_1^\#(r_k)\} = f(x) - L^0\{u_p^1(x)\} = b^1, \quad (30b)$$

$$\sum_{k=1}^{L} \lambda_k^n L^{n-1}\{u_n^\#(r_k)\} = L^{n-2}\{f(x)\} - L^{n-1}\{u_p^n(x)\} = b^n,$$
$$n=2,3,\ldots. \quad (30c)$$

In terms of the MRM, the successive process is truncated at some order $M$, namely, let

$$L^{M-1}\{u_p^M\} = 0. \quad (31)$$

The practical solution procedure is a reversal recursive process:

$$\lambda_k^M \to \lambda_k^{M-1} \to \cdots \to \lambda_k^0. \quad (32)$$

It is noted that due to

$$L^{n-1}\{v^n(r_k)\} = v^0(r_k), \quad (33)$$

The coefficient matrices of all successive equation are the same, i.e.

$$Q\lambda_k^n = b^n, \quad n=M, M-1,\ldots,1, \quad (34)$$

where matrix $Q$ is symmetric if operator $L\{\}$ is self-adjoint. Therefore, the LU decomposition algorithm is proper for this task. Finally, the particular solution $u_p^0$ is given by

$$u_p^0 = \sum_{n=1}^{M} \sum_{k=1}^{L} \lambda_k^n u_n^\#(r_k). \quad (35)$$

Substituting $u_p^0$ into Eq. (30a) yields

$$A\lambda_k^0 = b^0. \quad (36)$$

After getting the coefficient $\lambda_k^0$, we can calculate the solution at any interior or boundary node by

$$u(x_i) = \sum_{n=0}^{M} \sum_{k=1}^{L} \lambda_k^n u_n^\#(r_{ik}). \quad (37)$$

To differentiate the above solution procedure from the other boundary-type discretization schemes and manifest its meshfree merit, the approach is named as the boundary particle method (BPM). We can further divide the BPM as the BPM-1 due to singular fundamental solution and the BPM-2 due to nonsingular general solution, which respectively correspond to the MFS and BKM. It is noted that the BPM with $M=1$ degenerates into the BKM without using the inner nodes. The BPM-1 using the singular fundamental solution requires the fictitious auxiliary boundary outside physical domain just like the MFS. In addition, a mixed BPM and BKM (or MFS) strategy is also interesting, that is, instead of truncated Eq. (31), we apply the dual reciprocity and RBF to evaluate higher–order particular solution $u_p^M$ through equation

$$L^{M-1}\{u_p^M\} = L^{M-2}\{f(x)\}. \qquad (38)$$

This mixed method may require interior nodes.

It is noted that the only difference between the BKM (MFS) and BPM lies in how to evaluate the particular solution. The former applies the dual reciprocity principle, while the latter employs the multiple reciprocity principle. The advantage of the BPM over the BKM is that it dose not require interior nodes which may be especially attractive in some problems such as moving boundary, inverse problems, and exterior Helmholtz problems. However, the BPM is more mathematically complicated due to the use of higher-order fundamental or general solutions. In addition, the BPM may cost more computing effort than the BKM in nonlinear and time-dependent cases. It is expected that like the MRM [15,19], the truncated order $M$ in the BPM may not be large (usually two or three orders) in a variety of practical uses.

Noting the actual equivalence between BPM formulation (30a) using nonsingular general solution and BKM formulations (12) and (13), it is rather straightforward to derive the symmetric BPM scheme by replacing Eq. (30a) by Eqs. (16) and (17), which is based on the Hermite RBF interpolation (15). Furthermore, one can easily develop the symmetric direct BPM methodology by replacing Eq. (36) by Eq. (24).

It is stressed that both the BKM and BPM circumvent the troublesome singular integral inherent in the BEM. Therefore, they are a very competitive alternative to the latter for practical engineering computations.

## 3. New domain-type RBF schemes

### 3.1. Modified Kansa method based on the Green integral

Kansa [7] introduced the first domain-type RBF schemes which is now called as the Kansa's method. Despite great effort, the rigorous mathematical proof of the solvability of the Kansa's method is still missing [20]. One drawback in the Kansa's method is that the mixed boundary conditions may destroy the symmetric interpolation matrix. Fasshauer [8] presented the symmetric Hermite RBF collocation scheme which produces the symmetric matrix irrespective of the governing and boundary condition equations. In addition, the symmetricity also mathematically guarantees the solvability of Fasshauer's Hermite RBF scheme with some conditions.

One common problem in the Kansa's method and Fsshauer's symmetric Hermite method is that the numerical solutions at nodes adjacent to boundary are generally relatively much less accurate (by one to two orders) than those in the domain far from the boundary. Fedoseye et al. [21] proposed an improved Kansa-MQ scheme to effectively remove this shortcoming. The strategy is named as the PDE collocation on the boundary (PDECB). However, it is noted that the PDECB requires an additional set of nodes (inside or outside of the domain) adjacent to the boundary. Like the fictitious boundary in the MFS, the arbitrary placing of these additional nodes may give rise to some troublesome issues. In addition, the Kansa method with the PDECB lacks the explicit theoretical endorsement. In fact, the strategy similar to the PDECB has been independently proposed by Zhang et al. [22], which collocates both governing and boundary equations on the same boundary nodes. However, the method given in [22] is unsymmetrical and still lacks explicit theoretical basis.

By using the Green second identity, the present author presents a modified Kansa's method (MKM), which eliminates the above weaknesses in the Kansa's method with the PDECB [21] and the Hermite-type method given by Zhang et al [22]. We construct symmetric MKM via Fasshauer's Hermite interpolation. To better illustrate the idea behind the MKM, let us still consider cases described by Eqs. (1) and (2a,b), its Green integral solution is given by

$$u(x) = \int_\Omega f(z) u^*(x,z) d\Omega(z) + \int_\Gamma \left\{ u \frac{\partial u^*(x,z)}{\partial n(z)} - \frac{\partial u}{\partial n(z)} u^*(x,z) \right\} d\Gamma(z), \quad (39)$$

where $u^*$ is the fundamental solution of differential operator $L\{\}$. $z$ denotes source point. It is noted that the first and second terms of Eq. (39) are respectively equivalent to the particular and general solutions of Eq. (3). If a numerical integral scheme is used to analogize Eq. (39), we have

$$u(x) \cong \sum_{k=1}^{N+L} \omega(x, x_k) f(x_k) u^* + \sum_{k=N+1}^{N+L} Q(x, x_k) \left[ u \frac{\partial u^*}{\partial n} - \frac{\partial u}{\partial n} u^* \right], \quad (40)$$

where $\omega(x, x_j)$ and $Q(x, x_j)$ are the integration weight functions dependent on the integral schemes. The present author [23] presented the kernel RBF-creating strategy via the above formula (40). By analogy with the Fasshauer's Hermite scheme, we can construct the following RBF interpolation formula

$$u(x) = \sum_{k=1}^{L+N} \alpha_k L^*\{\varphi(r_k)\} + \sum_{j=1}^{L_d} \beta_j \varphi(r_j) + \sum_{j=L_d+1}^{L} \beta_j \left( -\frac{\partial \varphi(r_j)}{\partial n} \right) \quad (41)$$

where the boundary nodes are interpolated twice, and $L^*\{\}$ reverses some signs of odd-order derivatives in $L\{\}$ if the latter is not self-adjoint. Eliminating operator $L^*\{\}$ in the first term and negative sign in the last term, we have the representation of the RBF approximation given by Zhang et al [22]

$$u(x) = \sum_{k=1}^{L+N} \alpha_k \varphi(r_k) + \sum_{j=1}^{L_d} \beta_j \varphi(r_j) + \sum_{j=L_d+1}^{L} \beta_j \frac{\partial \varphi(r_j)}{\partial n}, \quad (42)$$

It is stressed here that the formulas (41) and (42) are different from the standard Kansa's RBF expansion

$$u(x) = \sum_{k=1}^{L+N} \gamma_k \varphi(r_k) \quad (43)$$

in that the formers interpolate two times in all boundary nodes. Also the RBF collocation scheme based on the formula (41) differs from the Kansa's method with the PDECB in that the present strategy does not require additional auxiliary nodes at all and is derived naturally from the Green second identity. Therefore, the ambiguities in the PDECB due to an arbitrary set of additional nodes and unclear theoretical backgrounds are eliminated at all. Contrast to the method given in Zhang et al [22], the MKM is symmetric and clearly based on the Green integral expansion.

Collocating Eqs. (1) and (2a,b) via approximate representation (41), we have

$$A \begin{Bmatrix} \beta \\ \alpha \end{Bmatrix} = \begin{Bmatrix} B\{x_i\} \\ f(x_l) \end{Bmatrix}, \quad (44)$$

where $B\{x_i\}$ represents the boundary conditions and interpolation matrix $A$ is

$$A = \begin{bmatrix} \varphi & -\frac{\partial \varphi}{\partial n} & L^*\{\varphi\} \\ \frac{\partial \varphi}{\partial n} & -\frac{\partial^2 \varphi}{\partial n^2} & \frac{\partial L^*\{\varphi\}}{\partial n} \\ L\{\varphi\} & -\frac{\partial L\{\varphi\}}{\partial n} & LL^*\{\varphi\} \end{bmatrix}. \quad (45)$$

Note that the source and response nodes at each pair of diagonally-symmetric block are reciprocal. Therefore, the matrix $A$ will be symmetric no matter what kinds of boundary equations are involved. In particular, we need to stress that matrix $A$ holds symmetric property with radically symmetric RBF $\varphi$ even if operator $L\{\}$ is not self-adjoint. In this case, we need to create the corresponding operator $L^*\{\}$ with a different sign before the odd-order derivative in operator $L\{\}$. For instance, consider the convection-diffusion operator

$$L\{u\} = D\nabla^2 u - v \bullet \nabla u - ku, \quad (46)$$

where $v$ denotes velocity vector; $D$ is the diffusivity coefficient; $k$ represents the reaction coefficient; and $x$ means multi-dimensional independent variable. In order to have symmetric RBF interpolation matrix (45) for this non-self-adjoint operator, we need to create the corresponding $L^*\{\}$

$$L^*\{u\} = D\nabla^2 u + v \bullet \nabla u - ku. \quad (47)$$

In addition, even for varying parameter problems, say varying parameter Helmholtz equation

$$\nabla^2\{u\} + S(x)u = f(x), \quad (48)$$

where $S(x)$ is assumed positive, the MKM still keep symmetric matrix by modifying radial basis function as

$$\hat{\varphi}(r_{ik}) = \sqrt{S(x_i)S(x_k)}\varphi(r_{ik}). \quad (49)$$

It is noted that the MKM interpolation matrix is ($N+2L$) dimensions for the second order PDE's and ($N+3L$) dimensions for the fourth order PDE's.

The present author [6] proved that if the geometry of interested problems is symmetric, the RBF interpolation matrix has the centrosymmetric structure, which can be decomposed into two smaller size sub-matrices in the calculation of the determinant, inverse, and eigenvalue and eigenvectors. Therefore, the MKM matrix has symmetric centrosymmetric structure if the domain possesses the symmetric geometry such as rectangle and ellipse. Such matrix structure makes the MKM conserve many important physical properties of the system [24].

The difference between the present methodology and the Fasshauer's Hermite method lies in that we collocates the governing and boundary equations separately at each boundary node. We name the present RBF collocation scheme as the modified Kansa's method to differentiate it from the other domain-type schemes. Just like the previous boundary-type BKM and BPM, the domain-type MKM is also based on a simple fact that the solutions of a PDE system can be seen as a sum of particular and homogeneous solutions.

Another key issue concerning the accuracy and efficiency of the RBF collocation method is how to create basis function $\varphi$ in expansion series (41). We strongly recommend to employ the nonsingular higher-order fundamental solution or general solution relating to operator $L\{\}$ as the RBF in the domain-type schemes. For example, the higher-order fundamental solutions of 2D Laplace operator [15] are given by

$$u_j^* = \frac{1}{2\pi} r^{2j}(A_j \ln r - B_j), \quad (50)$$

where

$$A_{j+1} = \frac{A_j}{4(j+1)^2}, \quad (51a)$$

$$u_j^* = \frac{1}{4(j+1)^2}\left(\frac{A_j}{j+1} + B_j\right). \quad (51b)$$

$j=1$ represents the original Laplace operator. For Helmholtz operators, we have either higher-order fundamental solutions or nonsingular general solutions.

It is also expected that the MKM is less accurate than the BKM due to the fact that the latter employs the analytical nonsingular general solution to approximate homogeneous solution. However, the MKM can produce the symmetric interpolation matrix no matter whether the right-hand inhomogeneous term in Eq. (1) includes the unknown dependent variable $u$ and whether operator $L\{\}$ is self-adjoint or not. In addition, the MKM is also relatively easier to program.

### 3.2. Spline RBF interpolation and local MKM

Compared with the functional integral form of the FEM including boundary conditions, one can note that the MKM applies the boundary conditions in a similar fashion. In terms of the recently-developed Green element method [25], it is promising to further localize the MKM, which is the topic of this subsection.

The basic strategy localizing the MKM is to apply spline RBF interpolation and domain decomposition, which is also interpreted as influential domain, namely, the solution at one node within a specified subregion is approximated only by those at the nodes of the same subregion. It is stressed here that this subregion can by no way been understood as the elements or grids in the FEM and FDM. Therefore, this method is still truly meshless. In addition, the strategy does not involve the overlapping across different subregions. We call this localizing MKM as the finite knot method (FKM).

For example, if we require $C_0$ continuity across different internal subregions, we have the RBF interpolation expression

$$u(x) = \sum_{k=1}^{L_\zeta + N_\zeta} \alpha_k L^*\{\varphi(r_k)\} + \sum_{j=1}^{L_\zeta} \beta_j \varphi(r_j), \quad (52)$$

where $\zeta$ is index of subregion, and $N$ and $L$ with subscript $\zeta$ denotes the numbers of internal and boundary nodes of the corresponding subregion. If some part of boundary of a subregion involves the exterior physical boundary, we should use a different representation like Eq. (41) considering different types of boundary conditions.

For $C_1$ continuity at internal boundary nodes of different 2D subregions, we have

$$u(x) = \sum_{k=1}^{L_\zeta + N_\zeta} \alpha_k L^*\{\varphi(r_k)\} + \sum_{j=1}^{L_\zeta} \beta_j \varphi(r_j) + \sum_{m=1}^{L_\zeta} \gamma_m \left(-\frac{\partial \varphi(r_m)}{\partial x}\right) + \sum_{n=1}^{L_\zeta} \eta_n \left(-\frac{\partial \varphi(r_n)}{\partial y}\right) \quad (53)$$

Similarly we can construct the RBF spline with higher-order continuity. Based on the Green integral, the FKM is a Hermite spline local RBF interpolation scheme. By using this piecewise RBF interpolation, we can have sparse banded system matrix of the FKM discretizing PDEs.

It is worth pointing out that unlike the FEM, the FKM holds the symmetric matrix with radically symmetric radial basis function even if the differential operator lacks the self-adjoint property. In addition, the FKM eases severe ill-conditioning of the global RBF schemes [26]. Therefore, the FKM is a very competing technique to handle a broader range of large-size problems.

By applying local symmetric geometry partitioning, the factorization merit of symmetric centrosymmetric matrix can lead to further significantly reducing the computational effort, preserving the physical features of real system and improving computation stability in the FKM computing large-size problems [24].

### 3.3. Direct MKM and FKM schemes

According to interpolation expression (41), the approximation solutions to $u$ at all nodes of domain and boundary and to the Neumann boundary conditions at all boundary nodes can be expressed as

$$B\begin{Bmatrix}\beta\\\alpha\end{Bmatrix} = \begin{Bmatrix}u\\q\end{Bmatrix}, \quad (54)$$

where

$$B = \begin{bmatrix} \varphi & -\dfrac{\partial \varphi}{\partial n} & L^*\{\varphi\} \\ \dfrac{\partial \varphi}{\partial n} & -\dfrac{\partial^2 \varphi}{\partial n^2} & \dfrac{\partial L^*\{\varphi\}}{\partial n} \end{bmatrix}, \quad (55)$$

$$q = \frac{\partial u}{\partial n}. \quad (56)$$

$q$ is the Neumann boundary conditions. Note that the fist row block in matrix $B$ is different from that of matrix $A$ of Eq. (44) in that it includes the representation of function value $u$ at both boundary and domain nodes.

In terms of Eqs. (44) and (54), we have direct RBF approximation

$$\begin{Bmatrix}u\\q\end{Bmatrix} = BA^{-1}\begin{Bmatrix}B\{x_i\}\\f(x_l)\end{Bmatrix}. \quad (57)$$

It is seen from Eq. (57) that no indirect expansion coefficients are used there. This direct MKM formulation is advantageous for some practical problems, especially when encountering discontinuous solution such as shock. It is straightforward to derive the direct FKM formulation.

## 3.4. Iterative RBFs for nonlinear and varying-parameter problems

For varying-parameter and nonlinear problems, constructing efficient RBF becomes a daunting task. A feasible strategy is to iterative alternation of the RBFs at different nodes through the whole iterative solution process, which is called as iterative RBFs. For example, let us consider nonlinear steady 2D Burger equation

$$\nabla^2 u - u_x u = f(x), \qquad (58)$$

the effective RBF should be the nonsingular higher-order fundamental or general solution of convection-diffusion or modified Helmholtz operator. In terms of bi-modified Helmholtz operator, we have the RBF

$$u(r,x) = rI_1\left(r\sqrt{|u_x|}\right)/2|u_x|^3 + I_0\left(r\sqrt{|u_x|}\right)/4u_x^4, \qquad (59)$$

where $I_0$ and $I_1$ are respectively the modified Bessel functions of the first kind of the zero and first orders. The approximate general solution for bi-convection diffusion operator can be written as

$$u(r,x) = e^{-u(x-x_k)/2} rI_1\left(\frac{|u|}{\sqrt{2}} r\right), \qquad (60)$$

The iterative RBF means that we need to alternate progressively the varying parameter $u_x$ or $u$ in RBF (59) or (60) by their numerical solutions of the previous iterative step.

## 4. Some remarks

The BKM and BPM are global interpolation techniques and always produce full matrix. In contrast, the BEM is not a global scheme with lower order of accuracy but still encounters the full matrix. Therefore, the BKM and BPM are theoretically more attractive than the BEM. Nonetheless, a sparse banded interpolation matrix is still desirable to ease severe ill-conditioning of large dense matrix, especially when an iterative method is used in the solution of linear or nonlinear discretization equations. A feasible sparse strategy is to use the spline interpolation as in the Green element method [25]. Namely, we apply the BKM or BPM in sub-regions and assemble all local approximations into the resultant global matrix by spline continuity conditions. It is noted that the spline BKM and BPM are not element-based techniques and truly meshless.

Beatson's rapid RBF algorithms [27] based on the multipole and moment methods were claimed capable of solving huge-size RBF interpolation systems of millions nodes on a modest desktop PC. The rapid algorithms can be used to significantly reduce the computational effort of the BKM, BPM, and MKM and enable these methods for practically significant problems. The popular wavelets-based rapid solution [28] of algebraic equations can also be employed to handle dense RBF interpolation matrix.

In addition, the underlying relationship between the RBF and wavelet recently unveiled by Chen [6,29] show that the RBF has inherent capability handling multiscale problems. The special matrix product techniques for nonlinear collocation formulation developed by Chen et al. [18] is also very efficient in the solution of RBF discretization equations of nonlinear PDEs.

In this study, we only discuss the collocation RBF schemes. Combining the presented RBF boundary-type and domain-type interpolation expressions with the Galerkin or least squares discretization are also interesting research topics.

The numerical experiments applying the presented BKM, BPM, MKM and FKM will be reported in subsequent papers. It is stressed that all these novel RBF collocation schemes proposed here are direct, symmetric, meshless, and integration-free. The spline-based schemes will produce desirable symmetric sparse banded interpolation matrix. The compactly-supported RBF is important to the spline RBF discretization.

## Appendix

We here present a Hermite RBF interpolation strategy to eliminate so-called edge effect in geometry reconstruction and image processing [23,30], that is, the interpolation solution at the nodes adjacent to edge are relatively less accurate then that at the nodes far away from edge. In terms of the MKM, we have RBF interpolation representation

$$u(x) = \sum_{k=1}^{L+N} \alpha_k \varphi(r_k) + \sum_{j=1}^{L} \beta_j \left( -\frac{\partial \varphi(r_j)}{\partial n} \right) \quad (A1)$$

where $L$ and $N$ are respectively the numbers of boundary and interior nodes, and the second term denotes the Numann boundary data, if it is available. Of course, we can replace Neumann data by the other type boundary data in (A1). The interpolation equation can be written as

$$\begin{bmatrix} \varphi & -\frac{\partial \varphi}{\partial n} \\ \frac{\partial \varphi}{\partial n} & -\frac{\partial^2 \varphi}{\partial n^2} \end{bmatrix} \begin{Bmatrix} \alpha \\ \beta \end{Bmatrix} = \begin{Bmatrix} u \\ \frac{\partial u}{\partial n} \end{Bmatrix}. \quad (A2)$$

Note that the above interpolation matrix is still symmetric. The essence of this methodology is to enforce additional boundary restraints, i.e., Neumann data in (A1), on the geometric modeling or imaging. The present strategy requires more boundary information than the existing methods. Beatson's rapid algorithm [27] can be applied to the interpolation equation (A2).

**References**

[1] Belytschko, T., Lu, YY. and Gu, L. Element-free Galerkin methods, *Int. J. Numer. Methods Eng*. 37, 1994, 229-256.

[2] De, P., Mendocan, T.R., de Barcellos, C.S. and Duarte, A., Investigations on the hp-Cloud Method by solving Timoshenko beam problems. *Comput. Mech.* 25, 2000, 286-295.

[3] Li, S. and Liu, W.K., Meshfree and particle methods and their applications, *Applied Mechanics Review*, (to appear) 2001.

[4] Buhmann, M.D., Radial basis functions, *Acta Numerica*, 2000, 1-38.

[5] Chen, W. and Tanaka, M., New Insights into Boundary-only and Domain-type RBF Methods, *Int. J. Nonlinear Sci. & Numer. Simulation*, 1(3), 2000, 145-151.

[6] Chen, W. and Tanaka, M., Relationship between boundary integral equation and radial basis function, *The 52th Symposium of JSCME on BEM*, (Tanaka, M. ed.), Tokyo, 2000.

[7] Kansa, E.J., Multiquadrics: A scattered data approximation scheme with applications to computational fluid-dynamics, *Comput. Math. Appl.* 19, 1990, 147-161.

[8] Fasshauer, G.E., Solving partial differential equations by collocation with radial basis functions, Proceedings of Chamonix, 1996, (Mehaute, A. and Rabut, C. and Schumaker, L. ed.), 1-8.

[9] Golberg, M.A. and Chen, C.S., The method of fundamental solutions for potential, Helmholtz and diffusion problems. In *Boundary Integral Methods - Numerical and Mathematical Aspects,* (Golberg, M.A. ed.), pp. 103-176, Comput. Mech. Publ., 1998.

[10] Nardini, D. and Brebbia, C.A., A new approach to free vibration analysis using boundary elements. *Appl. Math. Modeling*, 7, 1983, 157-162.

[11] Franke, C. and Schaback, R., Solving partial differential equation by collocation using radial basis function. *Appl. Math. Comput.* 93, 1998, 73-82.

[12] Kitagawa, T., Asymptotic stability of the fundamental solution method. *J. Comput. Appl. Math.*, 38, 1991, 263-269.

[13] Hon, Y. C. and Chen, W., Boundary knot method for 2D and 3D Helmholtz and convection-diffusion problems with complicated geometry, *Int. J. Numer. Methods Engng.*, (submitted), 2001.

[14] Golberg, M.A., Chen, C.S., Bowman, H. and Power, H., Some comments on the use of radial basis functions in the dual reciprocity method. *Comput. Mech* 21, 1998, 141-148.

[15] Nowak, A.J. and Neves, A.C. (eds), *The Multiple Reciprocity Boundary Element Method*, Comput. Mech. Publ., Southampton, U.K., 1994.

[16] Partridge, P.W., Brebbia, C.A. and Wrobel, L.W., *The Dual Reciprocity Boundary Element Method*, Comput. Mech. Publ., Southampton, UK, 1992.

[17] Kane, J.H., *Boundary Element Analysis in Engineering Continuum Mechanics*, Prentice Hall, 1993.

[18] Chen, W., Shu, C., He, W. and T. Zhong, The

DQ solution of geometrically nonlinear bending of orthotropic rectangular plates by using Hadamard and SJT product, *Computers & Structures,* 74(1), 2000, 65-74.
[19] Nowak, A.J., The multiple reciprocity method of solving heat conduction problems, Advances in Boundary Elements Method Vol.2 Field and Fluid Flow Solutions, (Brebbia, C.A. and Connor, J.J. eds.), pp. 81-95, Proc. Of the eleventh Int. Confer. Boundary Element Methods, Comput. Mech. Publ. & Springer-Verlag, 1989.
[20] Schaback, R. and Hon, Y.C., On unsymmetric collocation by radial basis functions, *J. Appl. Math. Comp.* 119, 2001, 177-186.
[21] Fedoseyev, A.I., Friedman, M.J. and Kansa, E.J., Improved multiquadratic method for elliptic partial differential equations via PDE collocation on the boundary. *Comput. Math. Appl.*, (in press), 2000.
[22] Zhang, X. Song, K.Z., Lu, M.W. and Liu, X., Meshless methods based on collocation with radial basis functions, *Comput. Mech.* 26, 2000, 333-343.
[23] Chen, W. and He, J., A study on radial basis function and quasi-Monte Carlo methods, *Int. J. Nonlinear Sci. & Numer. Simulation,* 1(4), 2000, 337-342.
[24] Cantoni, A. and Butler, P., Eigenvalues and eigenvectors of symmetric centrosymmetric matrices. *Linear Algebra Appl.* 13, 1976, 275-288.
[25] Taigbenu, AE, Enhancing the accuracy of the solution to unsaturated flow by a Hermitian Green element model, Adv. Eng. Software, 29(2), 1998, 113-118.
[26] Kansa, E.J. and Hon, Y.C., Circumventing the ill-conditioning problem with multiquadric radial basis functions: applications to elliptic partial differential equations. *Comput. Math. Appls.* 39, 2000, 123-137.
[27] Beatson, R.K., Cherrie, J.B. and Ragozin, D.L., Fast evaluation of radial basis functions: methods for four-dimensional polyharmonic splines, *SIAM J. Math. Anal.*, 32(6), 2001, 1272-1310.
[28] Beylkin, G., Coifman, R. and Rokhlin, V., Fast wavelet transforms and numerical algorithms: I, *Comm. Pure Appl. Math.*, 44, 1991, 141-183.
[29] Chen, W., Orthonormal RBF wavelet and ridgelet-like series and transforms for high-dimensional problems, *Int. J. Nonlinear Sci. Numer. Simulation*, 2(2), 2001, 163-168.
[30] Powell, M., Recent advances in Cambridge on radial basis functions. *The 2nd Int. Dortmunt Meeting on Approximation Theory*, 1998.